\documentclass[11pt,bezier]{article}
\usepackage{amsmath}
\usepackage{amsfonts,amsthm,amssymb}
\usepackage{amsfonts}
\usepackage{graphics}
\usepackage{amssymb}
\newtheorem{lem}{Lemma}[section]
\newtheorem{thm}[lem]{Theorem}

\newtheorem{conj}{Conjecture}

\theoremstyle{definition}

\begin{document}
\title{Nonseparating trees in 2-connected graphs and oriented trees in strongly connected digraphs
\footnote{The research is supported by NSFC (Nos.11401510, 11531011) and NSFXJ(No.2015KL019).}}
\author{Yingzhi Tian$^{a}$ \footnote{Corresponding author. E-mail: tianyzhxj@163.com (Y.Tian),
hjlai@math.wvu.edu (H. Lai), xuliqiong@jmu.edu.cn (L. Xu), mjx@xju.edu.cn (J.Meng).}, Hong-Jian Lai$^{b}$, Liqiong Xu$^{c}$, Jixiang Meng$^{a}$ \\
{\small $^{a}$College of Mathematics and System Sciences, Xinjiang
University, Urumqi, Xinjiang 830046, PR China}\\
{\small $^{b}$Department of Mathematics, West Virginia University,
Morgantown, WV 26506, USA}\\
{\small $^{c}$School of Science, Jimei University, Xiamen, Fujian 361021, PR China}}

\date{}

\maketitle

\noindent{\bf Abstract } Mader [J. Graph Theory 65 (2010) 61-69] conjectured that for every positive integer $k$ and every finite tree $T$ with order $m$,  every $k$-connected, finite graph $G$ with $\delta(G)\geq \lfloor\frac{3}{2}k\rfloor+m-1$ contains a subtree $T'$ isomorphic to $T$
such that $G-V(T')$ is $k$-connected. The conjecture has been verified for paths, trees when $k=1$, and stars or double-stars when $k=2$. In this paper we verify the conjecture for two classes of trees when $k=2$.

For digraphs, Mader [J. Graph Theory 69 (2012) 324-329] conjectured that every $k$-connected digraph $D$ with minimum semi-degree $\delta(D)=min\{\delta^+(D),\delta^-(D)\}\geq 2k+m-1$ for a positive integer $m$ has a dipath $P$ of order $m$  with $\kappa(D-V(P))\geq k$.  The conjecture has only been verified for the dipath with $m=1$, and the dipath with $m=2$ and $k=1$. In this paper, we prove that every strongly connected digraph with minimum semi-degree $\delta(D)=min\{\delta^+(D),\delta^-(D)\}\geq m+1$ contains an oriented tree $T$ isomorphic to some given oriented stars or double-stars with order $m$ such that $D-V(T)$ is still strongly connected.


\section{Introduction}

In this paper, all graphs (digraphs) are finite and without multiple edges (parallel arcs) and without loops. For graph-theoretical terminologies and notation not defined here, we follow \cite{Bondy}. A graph (digraph) is $k$-connected  means (strongly) $k$-vertex-connected. We use $\kappa(G)$ ($\kappa(D)$) to denote the connectivity of the graph $G$ (digraph $D$). The $order$ of a graph $G$ (digraph $D$) is the cardinality of its vertex set, denoted by $|G|$ ($|D|$).

In 1972, Chartrand, Kaugars, and Lick proved the following well-known result.

\begin{thm}{\cite{Chartrand}} Every $k$-connected graph $G$ of minimum degree
$\delta(G)\geq \lfloor\frac{3}{2}k\rfloor$ has a vertex $u$ with $\kappa(G-u)\geq k$.
\end{thm}

Fujita and Kawarabayashi proved in \cite{Fujita} that every $k$-connected graph $G$ with minimum degree at least $\lfloor\frac{3}{2}k\rfloor+2$ has an edge $e=uv$ such that $G-\{u,v\}$ is still $k$-connected. In the same paper, they stated the following conjecture.

\begin{conj}{\cite{Fujita}} For all positive integers $k, m$, there
is a (least) non-negative integer $f_k(m)$ such that every $k$-connected graph $G$ with $\delta(G)\geq \lfloor\frac{3}{2}k\rfloor-1+f_k(m)$ contains a connected subgraph $W$ of exact order $m$ such that $G-V(W)$ is still $k$-connected.
\end{conj}

The examples given in \cite{Fujita} showed that $f_k(m)$ must be at least $m$ for all positive integers $k, m$. In \cite{Mader2}, Mader confirmed Conjecture 1 and  proved that $f_k(m)=m$ holds for all $k, m$.

\begin{thm}{\cite{Mader2}}
Every $k$-connected graph $G$ with $\delta(G)\geq\lfloor\frac{3}{2}k\rfloor+m-1$
for positive integers $k, m$ contains a path $P$ of order $m$ such that $G-V(P)$ remains $k$-connected.
\end{thm}

Mader \cite{Mader2} further conjectured that Theorem 1.2 holds for all trees.

\begin{conj}{\cite{Mader2}} For every positive integer $k$ and every finite tree $T$, there is a least non-negative integer $t_k(T)$, such that every $k$-connected, finite graph $G$ with $\delta(G)\geq \lfloor\frac{3}{2}k\rfloor-1+t_k(T)$ contains a subgraph $T'\cong T$ with $\kappa(G-V(T'))\geq k$. Furthermore, $t_k(T)=|T|$ holds.
\end{conj}

In \cite{Mader3}, Mader showed that $t_k(T)$ exists. Actually, he showed that $t_k(T)\leq 2(k+m-1)^2+m-\lfloor\frac{3}{2}k\rfloor$. Theorem 1.2 implied that Conjecture 2 is true when $T$ is a path. Diwan and Tholiya \cite{Diwan} proved that Conjecture 2 holds when $k=1$.  In \cite{Tian1}, the authors verified that Conjecture 2 is true when $T$ is a star or double-star and $k=2$. In Section 3, we will verify Conjecture 2 for two classes of trees when $k=2$.

The minimum outdegree and the minimum indegree of a digraph $D$ are denoted by $\delta^+(D)$ and $\delta^-(D)$, respectively. The minimum semi-degree of $D$ is $\delta(D):=min\{\delta^+(D),\delta^-(D)\}$. The following result is a digraph analogue to Theorem 1.1.

\begin{thm}{\cite{Mader1}} Every $k$-connected digraph $D$ with minimum semi-degree
$\delta(D)=min\{\delta^+(D),\\ \delta^-(D)\} \geq 2k$ has a vertex $u$ with $\kappa(D-u)\geq k$.
\end{thm}

Considering the results for graphs and digraphs, Mader \cite{Mader3} suggested the following conjecture.

\begin{conj}{\cite{Mader3}} Every $k$-connected digraph $D$ with minimum semi-degree $\delta(D)=min\{\delta^+(D),\\ \delta^-(D)\} \geq 2k+m-1$ for a positive integer $m$ has a dipath $P$ of order $m$ with $\kappa(D-V(P))\geq k$.
\end{conj}

Mader remarked that one could conjecture also similar results for trees with special orientations, but he thought even a proof of Conjecture 3 very difficult. Conjecture 3 has only been verified for the dipath with $m=1$, and the dipath with $m=2$ and $k=1$. In Section 4, we will prove that every strongly connected digraph with minimum semi-degree $\delta(D)=min\{\delta^+(D),\delta^-(D)\}\geq m+1$ contains an oriented tree $T$ isomorphic to some given oriented stars or double-stars with order $m$ such that $D-V(T)$ is still strongly connected.

\section{Preliminaries}

Let $G$ be a graph with vertex set $V(G)$ and edge set $E(G)$. We write $u\in G$ for $u\in V(G)$. For a vertex $u\in G$, let $N_G(u)$ be the set of neighbors of $u$ in $G$ and $d_G(u)=|N_G(u)|$ be the $degree$ of $u$ in $G$.
For a vertex subset $U$ of a graph $G$, $G(U)$ denotes the subgraph induced by $U$ and $G-U$ is the subgraph induced by $V(G)-U$. The $neighborhood$ $N_G(U)$ of $U$ is the set of vertices in $V(G)-U$ which are adjacent to some vertex in $U$. If $U=\{u\}$, we use $G-u$ for $G-\{u\}$. If $H$ is a subgraph of $G$, we often use $H$ for $V(H)$. For example, $N_G(H)$, $H\cap G$ and $G(H)$ mean $N_G(V(H))$, $V(H)\cap V(G)$ and $G(V(H))$, respectively. If there is no confusion, we always delete the subscript, for example, $d(u)$ for $d_G(u)$, $N(u)$ for $N_G(u)$,  $N(U)$ for $N_G(U)$ and so on. For $H\subseteq G$, we define $\delta_G(H):=min_{x\in H}d_G(x)$, whereas $\delta(H)$ is the minimum degree of the graph $H$.  For $H_1,H_2\subseteq G$, $H_1\cup H_2$ is the subgraph of $G$ with vertex set $V(H_1)\cup V(H_2)$ and edge set $E(H_1)\cup E(H_2)$. For a set $S$, $K(S)$ denotes the complete graph on vertex set $S$.

A vertex set $S$ is a $separating$ $set$ of a connected graph $G$ if $G-S$ is disconnected, and $S$ is a $minimum$ $separating$ $set$ if $|S|=\kappa(G)$.
For a minimum separating set $S$ of $G$, we call the union $F$ of at least one, but not all components of $G-S$ a $fragment$ $F$ to $S$, and $\overline{F}:=G-(S\cup V(F))$ the $complementary$ $fragment$. An $end$ of $G$ is a fragment of $G$ which does not contain another fragment of $G$. An end of $G$ exists if and only if $G$ is not complete, and then, of course, there are at least two. The $completion$ of $S\subseteq V(G)$ in $G$, denoted by $G[S]$, is the graph $G\cup K(S)$.

Let $\mathcal{K}_k(m)$ denote the class of all pairs $(G,C)$, where $G$ is a $k$-connected graph with $|G|\geq k+1$, $C$ is a complete subgraph of $G$ with $|C|=k$ and with $\delta_G(G-V(C))\geq\lfloor\frac{3}{2}k\rfloor+m-1$. Let $\mathcal{K}_k^+(m)$ consist of all $(G,C)\in\mathcal{K}_k(m)$ with $\kappa(G)\geq k+1$.

In order to use induction to prove Theorem 1.2, Mader \cite{Mader2} proved the following result.

\begin{thm} (Mader \cite{Mader2})
For all $(G,C)\in \mathcal{K}_k^+(m)$ and $p\in G-V(C)$, there is a path $P\subseteq G-V(C)$ of order $m$ starting from $p$ such that $\kappa(G-V(P))\geq k$ holds.
\end{thm}

The following Theorem was stated in \cite{Mader3}. A proof of Theorem 2.2 was not given, but it follows from Theorem 2 in \cite{Mader2} in a similar way as that of Theorem 1 in \cite{Mader2}.

\begin{thm} (Mader \cite{Mader3})
Let $G$ be a $(k+1)$-connected graph $G$ with $\delta(G)\geq\lfloor\frac{3}{2}k\rfloor+m-1$ and let $p$ be a vertex of $G$. Then there is a path $P$ of order $m$ starting from $p$ such that $\kappa(G-V(P))\geq k$ holds.
\end{thm}

A $tree$ is a connected graph without cycles. A $star$ is a tree that has exact one vertex with degree greater than one. We call this vertex $u$ with degree greater than one the center-vertex of the star. A $double$-$star$ is a tree that has exact two vertices with degree greater than one. Those two vertices $u$ and $v$ with degree greater than one must be adjacent in a double-star. We call this edge $uv$ the center-edge of the double-star. The authors in \cite{Tian1} proved the following Theorem, which verified Conjecture 2 for stars and double-stars when $k=2$.

\begin{thm} \cite{Tian1}
Let $G$ be a 2-connected graph with minimum degree $\delta(G)\geq m+2$, where $m$ is a positive integer. Then

(i) $G$ contains a star $T'$ with order $m$ such that $G-V(T')$ is 2-connected;

(ii) for every double-star $T$ with order $m$, $G$ contains a double-star $T'$ isomorphic to $T$ such that $G-V(T')$ is 2-connected.
\end{thm}

As indicated in \cite{Mader2}, given a $(G,C)\in \mathcal{K}_k(m)$, we can obtain a graph $G'$ from $(G,C)$  such that $G'$ satisfies the condition $\delta(G')\geq\lfloor\frac{3}{2}k\rfloor+m-1$ by pasting together sufficiently many copies of $G$ at $C$. Then by applying Theorems 1.2 and 2.3 to $G'$, we obtain the following results.

\begin{thm}
(i) \cite{Mader2} Every $(G,C)\in \mathcal{K}_k(m)$ contains a path $P\subseteq G-V(C)$ of order $m$ such that $\kappa(G-V(P))\geq k$ holds.

(ii) \cite{Tian1} Every $(G,C)\in \mathcal{K}_2(m)$ contains a star $T'\subseteq G-V(C)$ of order $m$ such that $\kappa(G-V(T'))\geq 2$ holds;

(iii) \cite{Tian1} For every double-star $T$ with order $m$, every $(G,C)\in \mathcal{K}_2(m)$ contains a double-star $T'\subseteq G-V(C)$ isomorphic to $T$ such that $\kappa(G-V(T'))\geq 2$ holds.
\end{thm}

Let $D$ be a digraph with vertex set $V(D)$ and arc set $A(D)$. An arc $(u, v)$ is considered to be directed from $u$ to $v$. $u\in D$ means $u\in V(D)$.
For a vertex set $S\subseteq V(D)$, define $N_D^+(S):=\{v\in D-S:$ there is an $s\in S$ such that $(s,v)\in A(D)\}$; for a subdigraph $H\subseteq D$ and a vertex $u\in D$, we write $N_D^+(H)$ instead of $N_D^+(V(H))$ and $N_D^+(u)$ instead of $N_D^+(\{u\})$. For $u\in D$, $d_D^+(u)=|N_D^+(u)|$ denotes the outdegree of $u$. $\delta^+(D)$ denotes the minimum outdegree of $D$.
For the dual concepts, we use the notation $N_D^-$, $d_D^-$ and $\delta^-(D)$, respectively. The minimum semi-degree $\delta(D)$ of $D$ is $min\{\delta^+(D),\delta^-(D)\}$. The subdigraph of $D$ induced by $S\subseteq V(D)$ or $S\subseteq D$ is denoted by $D(S)$. We say $D$ is strongly connected if $\kappa(D)\geq1$.

\section{Connectivity keeping trees in 2-connected graphs}

The next lemma is widely used in studying connectivity of graphs.

\begin{lem} (Hamidoune \cite{Hamidoune})
Let $G$ be a $k$-connected graph and let $S$ be a separating set of $G$ with $|S|=k$. Then for every fragment $F$ of $G$ to $S$, $G[S]-V(\overline{F})$ is $k$-connected. Furthermore, if $F$ is an end of $G$ with $|F|\geq2$, then $G[S]-V(\overline{F})$ is $(k+1)$-connected.
\end{lem}

\noindent{\bf Proof.} Since a separating set of $G[S]-V(\overline{F})$ is one of $G$, too, Lemma 3.1 follows immediately. $\Box$

Before proving the main results in this section, we need one more lemma.

\begin{lem} (Mader \cite{Mader2})
Let $G$ be a $k$-connected graph and let $S$ be a separating set of $G$ with $|S|=k$. Then the following holds.

(i) Assume $\delta(G)\geq\lfloor\frac{3}{2}k\rfloor+m-1$ and let $F$ be a fragment of $G$ to $S$. If $W\subseteq G-(S\cup V(F))$ has order at most $m$ and $\kappa(G[S]-V(F\cup W))\geq k$ holds, then also $\kappa(G-V(W))\geq k$ holds.

(ii) Assume $(G,C)\in \mathcal{K}_k(m)$ and let $F$ be a fragment of $G$ to $S$ with $C\subseteq G(F\cup S)$. If $W\subseteq G-(S\cup V(F))$ has order at most $m$ and $\kappa(G[S]-V(F\cup W))\geq k$ holds, then also $\kappa(G-V(W))\geq k$ holds.
\end{lem}

\begin{center}
\scalebox{0.5}{\includegraphics{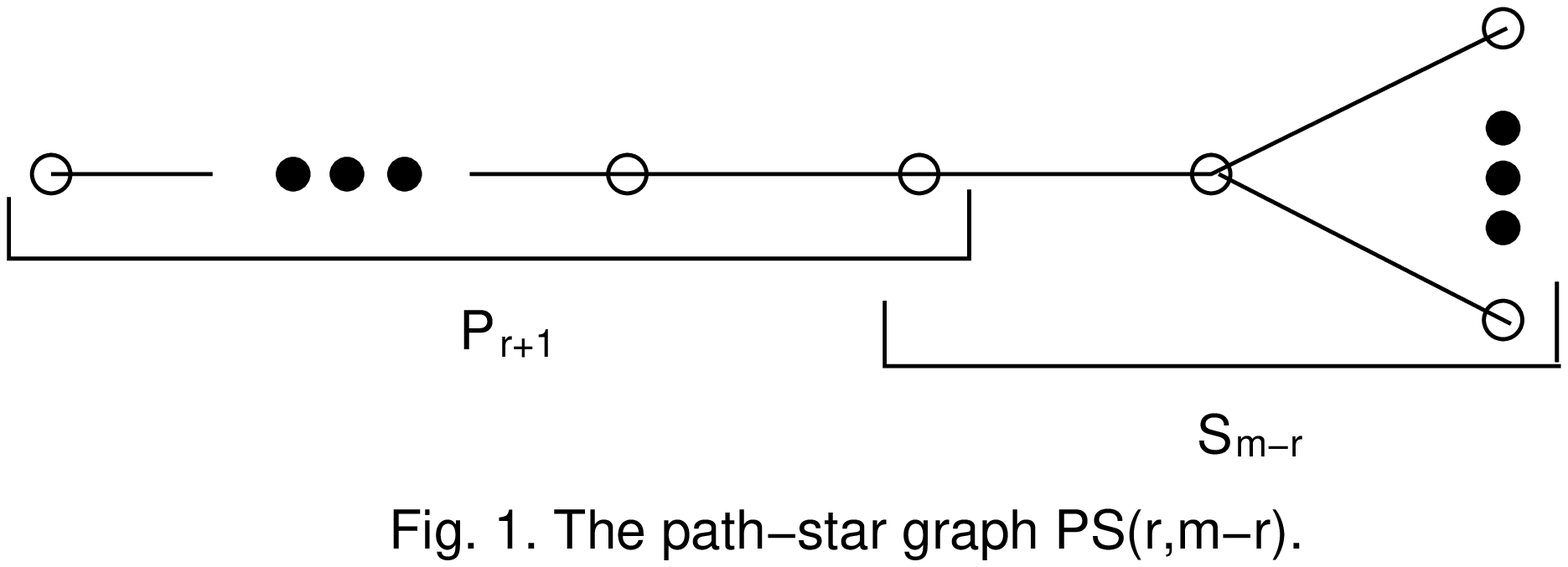}}
\end{center}

{\noindent {\bf Definition 1.}} The {\bf{path-star}} graph, $PS(r,m-r)$, is obtained from the disjoint union of  a path with order $r+1$ and a star with order $m-r$,  by identifying one end vertex of the path with one vertex degree one in the star. See Fig.1 for example.

\begin{thm}
Let $PS(r,m-r)$ be a path-star graph with order $m$, where $1\leq r\leq m-3$.
Then every 2-connected graph $G$ with minimum degree $\delta(G)\geq m+2$ contains a subgraph  $T$ isomorphic to $PS(r,m-r)$ such that $G-V(T)$ is 2-connected.
\end{thm}

\noindent{\bf Proof.} By Theorem 2.3(i), there is a star $T_1$ with order $m-r$ such that $G_1:=G-V(T_1)$ is 2-connected. Let $V(T_1)=\{u,v_1,\cdots,v_{m-r-1}\}$ and $E(T_1)=\{uv_i| i=1,\cdots,m-r-1\}$.
Since $\delta(G)\geq m+2$, we have $\delta(G_1)\geq m+2-(m-r)=r+2$ and $|N_G(v)\cap G_1|\geq m+2-(m-r-1)=r+3$ for each $v\in V(T_1)$.

\noindent{\bf Case 1.} $G_1$ is 3-connected.

Let $w_1$ be a neighbor of $v_1$ in $G_1$. By Theorem 2.2, there is a path $P_1$ of order $r$ in $G_1$ starting from $w_1$ such that $G_1-V(P_1)$ is 2-connected. Then the subgraph $T$ obtained from the union of $T_1$ and $P_1$ by adding an edge $v_1w_1$ satisfies $T\cong PS(r,m-r)$ and $G-V(T)$ is 2-connected.

\noindent{\bf Case 2.} $\kappa(G_1)=2$.

Then $G_1$ has an end $F$. Let $S=N_{G_1}(F)$ and $S=\{x,y\}$. We have $|F|\geq2$ by $\delta(G_1)\geq r+2$ and $|S|=2$. By Theorem 3.1, $G_1[S]-V(\overline{F})$ is 3-connected. By $\delta(G_1)\geq r+2$, we have $(G_1[S]-V(\overline{F}),S)\in \mathcal{K}_2^+(r)$.

\noindent{\bf Case 2.1.} $G_1$ has an end, say $F$, such that $N_G(T_1)\cap F\neq \O$.

If $v_i$ has a neighbor $w_i$ in $F$ for some $i\in\{1,\cdots,m-r-1\}$, then by Theorem 2.1, there is a path $P_1\subseteq G_1[S]-V(\overline{F})-S$ of order $r$ starting from $w_i$ such that $G_1[S]-V(\overline{F})-V(P_1)$ is 2-connected. Since $\delta(G_1)\geq r+2$ and $|V(P_1)|=r$, we have $G_1-V(P_1)$ is 2-connected by Lemma 3.2(i). Thus the subgraph $T$ obtained from the union of $T_1$ and $P_1$ by adding an edge $v_iw_i$ satisfies $T\cong PS(r,m-r)$ and $G-V(T)$ is 2-connected.

Assume $N_G(v_i)\cap F=\O$ for any $i\in\{1,\cdots,m-r-1\}$. By $N_G(T_1)\cap F\neq \O$, $u$ must have a neighbor, say $w$, in $F$. Let $T_2=T_1-v_{m-r-1}$ and $G_2=G-V(T_2)$. Then $F$ is still an end of $G_2$ and $\delta(G_2)\geq r+3$. Thus $(G_2[S]-V(\overline{F}),S)\in \mathcal{K}_2^+(r+1)$. By Theorem 2.1, there is a path $P_2\subseteq G_2[S]-V(\overline{F})-S$ of order $r+1$ starting from $w$ such that $G_2[S]-V(\overline{F})-V(P_2)$ is 2-connected. Since $\delta(G_2)\geq r+3$ and $|V(P_2)|=r+1$, we have $G_2-V(P_2)$ is 2-connected by Lemma 3.2(i). Thus the graph $T$ obtained from the union of $T_2$ and $P_2$ by adding an edge $uw$ satisfies $T\cong PS(r,m-r)$ and $G-V(T)$ is 2-connected.

\noindent{\bf Case 2.2.} Each end $F$ in $G_1$ satisfies $N_G(T_1)\cap F=\O$.

As above, assume $F$ is an end of $G_1$. Let $S=N_{G_1}(F)$ and $S=\{x,y\}$. Since $N_G(T_1)\cap F=\O$, we know $F$ is also an end of $G$. If we can find a subgraph $T\subseteq F$ such that $T\cong PS(r,m-r)$ and $\kappa(G[S]-V(\overline{F})-V(T))\geq2$, then, by applying Lemma 3.2(i) to $G$, we obtain $G-V(T)$ is 2-connected. Thus, in the following, we only need to prove that $G[S]-V(\overline{F})$ contains a subgraph $T'\subseteq G[S]-V(\overline{F})-S$ such that $T'\cong PS(r,m-r)$ and $\kappa(G[S]-V(\overline{F})-V(T'))\geq2$.

Let $G'=G[S]-V(\overline{F})$. By Lemma 3.1, $G'$ is 3-connected. Since $\delta(G)\geq m+2$, we have $(G',S)\in \mathcal{K}_2^+(m)$. By Theorem 2.4(ii), there is a star $T_1'\subseteq G'-S$ with order $m-r$ such that $G_1':=G'-V(T_1')$ is 2-connected. Let $V(T_1')=\{u',v_1',\cdots,v_{m-r-1}'\}$ and $E(T_1')=\{u'v_i'| i=1,\cdots,m-r-1\}$.
Since $\delta(G)\geq m+2$, we have $\delta_{G_1'}(G_1'-S)\geq m+2-(m-r)=r+2$ and $|N_{G'}(v')\cap G_1'|\geq m+2-(m-r-1)=r+3$ for each $v'\in V(T_1')$..

\noindent{\bf Case 2.2.1.} $G_1'$ is 3-connected.

Let $w_1'$ be a neighbor of $v_1'$ in $G_1'-S$ ($w_1'$ exists because $|N_{G'}(v_1')\cap G_1'|\geq r+3$ and $|S|=2$). By Theorem 2.1, there is a path $P_1'\subseteq G_1'-S$ of order $r$ starting from $w_1'$ such that $G_1'-V(P_1')$ is 2-connected. Then the graph $T'$ obtained from the union of $T_1'$ and $P_1'$ by adding an edge $v_1'w_1'$ satisfies $T'\cong PS(r,m-r)$ and $G'-V(T')$ is 2-connected.

\noindent{\bf Case 2.2.2.} $\kappa(G_1')=2$.

The proof of Case 2.2.2 is similar to Case 2.1. Nevertheless, we also outline the proof for completeness.

We can choose an end $F'$ of $G_1'$ such that $F'\cap S=\O$. Let $S'=N_{G_1'}(F')$ and $S'=\{x',y'\}$. We have $|F'|\geq2$ by $\delta_{G_1'}(G_1'-S)\geq r+2$ and $|S'|=2$. By Theorem 3.1, $G_1'[S']-V(\overline{F'})$ is 3-connected. By $\delta_{G_1'}(G_1'-S)\geq r+2$, we have $(G_1'[S']-V(\overline{F'}),S')\in \mathcal{K}_2^+(r)$.

If $v_i'$ has a neighbor $w_i'$ in $F'$ for some $i\in\{1,\cdots,m-r-1\}$, then by Theorem 2.1, there is a path $P_1'\subseteq G_1'[S']-V(\overline{F'})-S'$ of order $r$ starting from $w_i'$ such that $G_1'[S']-V(\overline{F'})-V(P_1')$ is 2-connected. Since $\delta_{G_1'}(G_1'-S)\geq r+2$ and $|V(P_1')|=r$, we have $G_1'-V(P_1')$ is 2-connected by Lemma 3.2(ii). Thus the graph $T'$ obtained from the union of $T_1'$ and $P_1'$ by adding an edge $v_i'w_i'$ satisfies $T'\cong PS(r,m-r)$ and $G'-V(T')$ is 2-connected.

Assume $N_{G'}(v_i)\cap F'=\O$ for any $i\in\{1,\cdots,m-r-1\}$. By $G'$ is 3-connected, $u'$ must have a neighbor, say $w'$, in $F'$. Let $T_2'=T_1'-u_{m-r-1}'$ and $G_2'=G'-V(T_2')$. Then $F'$ is still an end of $G_2'$ and $\delta_{G_2'}(G_2'-S)\geq r+3$. Thus $(G_2'[S']-V(\overline{F'}),S')\in \mathcal{K}_2^+(r+1)$. By Theorem 2.1, there is a path $P_2'\subseteq G_2'[S']-V(\overline{F'})-S'$ of order $r+1$ starting from $w'$ such that $G_2'[S']-V(\overline{F'})-V(P_2')$ is 2-connected. Since $\delta_{G_2'}(G_2'-S)\geq r+3$ and $|V(P_2')|=r+1$, we have $G_2'-V(P_2')$ is 2-connected by Lemma 3.2(ii). Thus the graph $T'$ obtained from the union of $T_2'$ and $P_2'$ by adding an edge $u'w'$ satisfies $T'\cong PS(r,m-r)$ and $G'-V(T')$ is 2-connected. $\Box$

{\noindent {\bf Definition 2.}} The {\bf{path-double-star}} graph, $PDS(r,m-r)$, is obtained from the disjoint union of a path with order $r+1$ and a double-star with order $m-r$,  by identifying one end vertex of the path with one vertex of degree one in the double-star. See Fig.2 for example.

 Specifically, we denote $PDS1(r,m-r)$ the path-double-star graph obtained from the disjoint union of a path with order $r+1$ and a double-star with order $m-r$,  by identifying one end vertex of the path with one pendant vertex which is adjacent to the vertex with maximum degree in the double-star. We denote $PDS2(r,m-r)$ the path-double-star graph obtained from the disjoint union of a path with order $r+1$ and a double-star with order $m-r$,  by identifying one end vertex of the path with one pendant vertex which is adjacent to the vertex with the second maximum degree in the double-star.

By replacing the star with a double-star in the proof of Theorem 3.3, we can obtain the proof of Theorem 3.4 by using almost the same arguments as the proof of Theorem 3.3. Besides, we use Theorem 2.3 (ii) instead of   Theorem 2.3 (i) and Theorem 2.4 (iii) instead of Theorem 2.4 (ii) in the proof of Theorem 3.4. So we omit the proof here.

\begin{thm}
For two integers $r,m$ with $1\leq r\leq m-4$, every 2-connected graph $G$ with minimum degree $\delta(G)\geq m+2$ contains a subgraph  $T$ isomorphic to $PDS1(r,m-r)$ or $PDS2(r,m-r)$ such that $G-V(T)$ is 2-connected.
\end{thm}

 \begin{center}
\scalebox{0.5}{\includegraphics{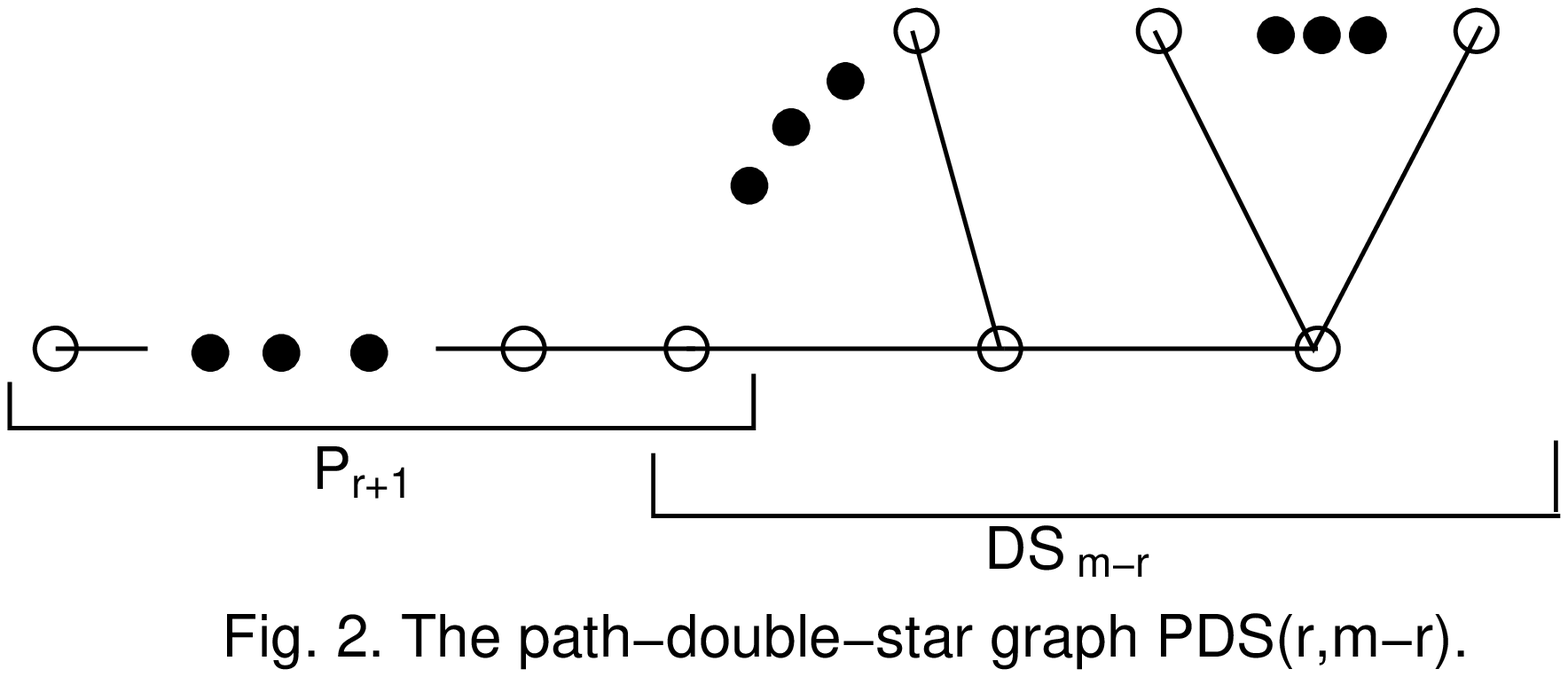}}
\end{center}

\section{Nonseparating oriented stars or double-stars in strongly connected digraphs}

{\noindent {\bf Definition 3.}} The {\bf{out-star}} $OS_m$ is the digraph  obtained from a star of order $m$ by orienting each edge of the star away from the center-vertex.  The {\bf{in-star}} $IS_m$ is the digraph  obtained from a star of order $m$ by orienting each edge of the star towards the center-vertex. See Fig.3 for examples.

\begin{center}
\scalebox{0.5}{\includegraphics{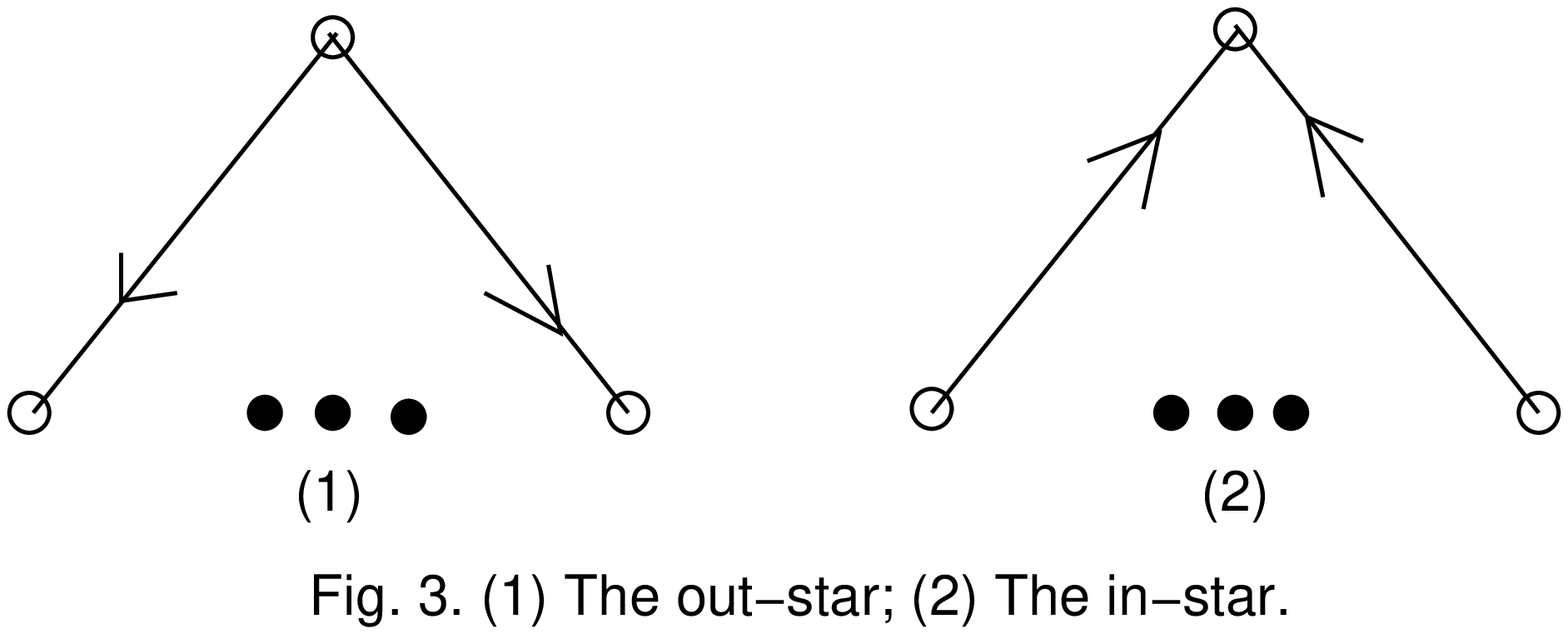}}
\end{center}

\begin{lem}
Let $D$ be a strongly connected digraph and $H$ be a subdigraph of $D$ with $|V(H)|<|V(D)|$. Then there is a dipath $P:=p_0p_1\cdots p_t(t\geq2)$ in $D$ such that $p_0,p_t\in V(H)$ and $p_1,\cdots,p_{t-1}\in V(D)-V(H)$, where $p_0$ and $p_t$ may be the same vertex.
\end{lem}

\noindent{\bf Proof.} Since $D$ is strongly connected, there is an arc, say $(p_0,p_1)$, from $V(H)$ to $V(D)-V(H)$. By $D$ is strongly connected, there is a dipath, say $P':=p_1\cdots p_t$, from $p_1$ to $V(H)$, where $p_1,\cdots,p_{t-1}\in V(D)-V(H)$  and $p_t\in V(H)$. Thus the dipath $P:=p_0p_1\cdots p_t$ is just a dipath we needed.  $\Box$

\begin{thm}
Let $D$ be a strongly connected digraph with minimum semi-degree $\delta(D)=min\{\delta^+(D), \delta^-(D)\}\geq m+1$. Then $D$ contains a subdigraph $T$ isomorphic to $OS_m$ or $IS_m$ such that $D-V(T)$ remains strongly connected.
\end{thm}

\noindent{\bf Proof.} Since $\delta(D)\geq m+1$, there are subdigraphs in $D$ isomorphic to $OS_m$ or $IS_m$.  Let $T$ be a subdigraph in $D$ isomorphic to $OS_m$ or $IS_m$. Let $D'=D-T$. If $D'$ is strongly connected, then we are done. Thus we assume that $D'$ is not strongly connected. We order all strong components of $D'$ as $C_1,\cdots,C_l$ such that there are no arcs from $C_j$ to $C_i$ for all $1\leq i<j\leq l$. Let $B$ be a maximum strong component of $D'$. We choose such a $T$ so that 

(1) $|B|$ is as large as possible.

Let $P:=p_0p_1\cdots p_t$ ($t\geq2$) be a shortest dipath in $D$ such that $p_0,p_t\in B$ and $p_1,\cdots,p_{t-1}\in V(D)-B$ (by Lemma 4.1). We consider three cases in the following.

\noindent{\bf Case 1.} $t=2$.

By $t=2$, we have $p_1\in V(T)$. If $B=C_1$, then for any vertex $c_l\in C_l$, we have $|N_D^+(c_l)\backslash (B\cup p_1)|\geq m+1-1=m$. Thus we can find an out-star $T'$ rooted at $c_l$ with order $m$ such that $V(T')\cap (V(B)\cup V(P))=\O$. But then $V(B)\cup V(P)$ is contained in a strong component of $D-V(T')$, contrary to (1). If $B\neq C_1$, then for any vertex $c_1\in C_1$, we have $|N_D^-(c_1)\setminus (B\cup p_1)|\geq m+1-1=m$. Thus we can find an in-star $T''$ rooted at $c_1$ with order $m$ such that $V(T'')\cap (V(B)\cup V(P))=\O$. But then $V(B)\cup V(P)$ is contained in a strong component of $D-V(T'')$, contrary to (1).

\noindent{\bf Case 2.} $t=3$.

By the choice of $P$, we have $N_D^+(p_1)\cap B=\O$. Then $|N_D^+(p_1)\setminus (B\cup P)|\geq m+1-1=m$. Let $q\in N_D^+(p_1)\setminus (B\cup P)$. If $N_D^+(q)\cap B=\O$, then $|N_D^+(q)\backslash (B\cup P)|\geq m+1-2=m-1$. Thus we can find an out-star $T'$ rooted at $q$ with order $m$ such that $V(T')\cap (V(B)\cup V(P))=\O$. But then $V(B)\cup V(P)$ is contained in a strong component of $D-V(T')$, contrary to (1). If $N_D^+(q)\cap B\neq\O$, then $N_D^-(q)\cap B=\O$ (for otherwise, we can find a dipath $P'$ shorter that $P$). Thus $|N_D^-(q)\backslash (B\cup P)|\geq m+1-2=m-1$, and we can find an in-star $T''$ rooted at $q$ with order $m$ such that $V(T'')\cap (V(B)\cup V(P))=\O$. But then $V(B)\cup V(P)$ is contained in a strong component of $D-V(T'')$, contrary to (1).

\noindent{\bf Case 3.} $t\geq4$.

By the choice of $P$, we have $N_D^+(p_1)\cap (B\cup P)=\{p_2\}$. Then $|N_D^+(p_1)\setminus (B\cup P)|\geq m+1-1=m$. Let $q_1,\cdots,q_m\in N_D^+(p_1)\setminus (B\cup P)$. By the choice of $P$, we have $N_D^+(q_j)\cap (B\cup P)\subseteq\{p_1,p_2,p_3\}$ for each $j\in\{1,\cdots,m\}$. Thus $|N_D^+(q_j)\setminus (B\cup P)|\geq m+1-3=m-2$ for each $j\in\{1,\cdots,m\}$.
If there is some $j\in\{1,\cdots,m\}$ such that $|N_D^+(q_j)\setminus (B\cup P)|\geq m-1$, then we can find an out-star $T'$ rooted at $q_j$ with order $m$ such that $V(T')\cap (V(B)\cup V(P))=\O$. But then $V(B)\cup V(P)$ is contained in a strong component of $D-V(T')$, contrary to (1). Thus we assume that $|N_D^+(q_j)\setminus (B\cup P)|=m-2$ and $N_D^+(q_j)\cap (B\cup P)=\{p_1,p_2,p_3\}$ for each $j\in\{1,\cdots,m\}$. But then there is an in-star $T''$ rooted at $p_2$ with vertex set $\{p_2,q_1,\cdots,q_{m-1}\}$ such that $V(T'')\cap (V(B)\cup V(P'))=\O$, where $P'$ is the dipath $p_0p_1q_mp_3\cdots p_t$. But then $V(B)\cup V(P')$ is contained in a strong component of $D-V(T'')$, contrary to (1).  $\Box$

{\noindent {\bf Definition 4.}} The {\bf{out-double-star}} $ODS(m;r,s)$ is the digraph obtained from the disjoint union of two out-stars (one is isomorphic to $OS_{r+1}$ rooted at $u$ and the other is isomorphic to $OS_{s+1}$ rooted at $v$, where $r+s=m-2$) by adding an arc from $u$ to $v$.
The {\bf{in-double-star}} $IDS(m;r,s)$ is the digraph obtained from the disjoint union of two in-stars (one is isomorphic to $IS_{r+1}$ rooted at $u$ and the other is isomorphic to $IS_{s+1}$ rooted at $v$, where $r+s=m-2$) by adding an arc from $u$ to $v$. The {\bf{out-in-double-star}} $OIDS(m;r,s)$ is the digraph obtained from the disjoint union of one out-star and one in-star (the out-star is isomorphic to $OS_{r+1}$ rooted at $u$ and the in-star is isomorphic to $IS_{s+1}$ rooted at $v$, where $r+s=m-2$) by adding an arc from $u$ to $v$. We often call the arc $(u,v)$ center-arc. See Fig.4 for examples.

\begin{center}
\scalebox{0.8}{\includegraphics{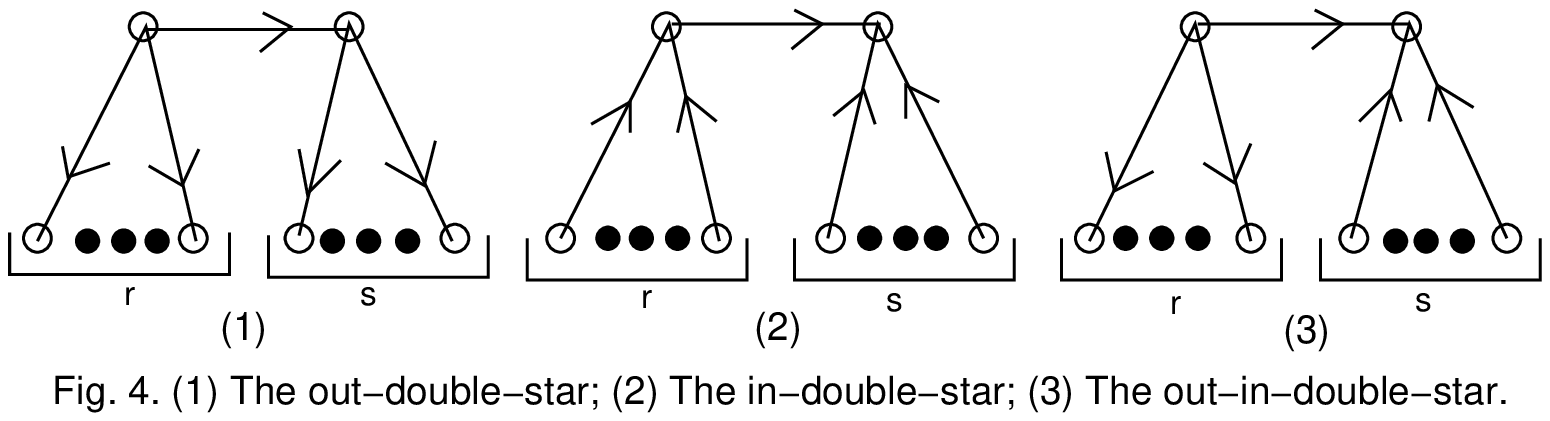}}
\end{center}

\begin{lem}
Let $D$ be a digraph, $T_1$ be an out-double-star with order $m$, $T_2$ be an in-double-star with order $m$, and $T_3$ be an out-in-double-star with order $m$.

(i) If there is an arc $a=(u,v)\in A(D)$ such that $|N_D^+(u)\setminus v|\geq m-3$, $|N_D^+(v)\setminus u|\geq m-3$ and $|N_D^+(u)\cup N_D^+(v)\setminus \{u,v\}|\geq m-2$, then there is an out-double-star $T\subseteq D$ isomorphic to $T_1$.

(ii) If there is an arc $a=(u,v)\in A(D)$ such that $|N_D^-(u)\setminus v|\geq m-3$, $|N_D^-(v)\setminus u|\geq m-3$ and $|N_D^-(u)\cup N_D^-(v)\setminus \{u,v\}|\geq m-2$, then there is an in-double-star $T\subseteq D$ isomorphic to $T_2$.

(iii) If there is an arc $a=(u,v)\in A(D)$ such that $|N_D^+(u)\setminus v|\geq m-3$, $|N_D^-(v)\setminus u|\geq m-3$ and $|N_D^+(u)\cup N_D^-(v)\setminus \{u,v\}|\geq m-2$, then there is an out-in-double-star $T\subseteq D$ isomorphic to $T_3$.
\end{lem}

\noindent{\bf Proof.} (i) By $T_1$ is an out-double-star, we have $m\geq4$. Assume the out-double-star $T_1$ has the center-arc $a'=(u',v')$, where $|N_{T_1}^+(u')\setminus v'|=r$, $|N_{T_1}^+(v')\setminus u'|=s$ and $|N_{T_1}^+(u')\cup N_{T_1}^+(v')\setminus \{u',v'\}|=m-2$ ($1\leq r,s\leq m-3$ and $r+s=m-2$). Since $|N_D^+(u)\setminus v|\geq m-3$, $|N_D^+(v)\setminus u|\geq m-3$ and $|N_D^+(u)\cup N_D^+(v)\setminus \{u,v\}|\geq m-2$, we can find an out-double-star $T\cong T_1$ in $D$ with center-arc $a=(u,v)$.

(ii) and (iii) follow similarly. $\Box$

While the main idea of the proof of Theorem 4.4 is similar to that of Theorem 4.2, the proof of Theorem 4.4 is somewhat more complicated with different details.

\begin{thm}
For any integers $m,r,s$ with $1\leq r,s\leq m-3$ and $r+s=m-2$, every
strongly connected digraph $D$ with minimum semi-degree $\delta(D)=min\{\delta^+(D), \delta^-(D)\}\geq m+1$ contains a subdigraph $T$ isomorphic to a member in ${\mathcal S} (m;r,s)=\{ODS(m;r,s), IDS(m;r,s), OIDS(m;r,s)\}$ such that $D-V(T)$ remains strongly connected.
\end{thm}

\noindent{\bf Proof.} Since $\delta(D)\geq m+1$, there is a subdigraph in $D$ isomorphic to a member in ${\mathcal S} (m;r,s)=\{ODS(m;r,s), IDS(m;r,s), OIDS(m;r,s)\}$ by Lemma 4.3.  Let $T$ be a subdigraph in $D$ isomorphic to a member in ${\mathcal S} (m;r,s)$. Let $D'=D-T$. If $D'$ is strongly connected, then we are done. Thus we assume that $D'$ is not strongly connected. We order all strong components of $D'$ as $C_1,\cdots,C_l$ such that there are no arcs from $C_j$ to $C_i$ for all $1\leq i<j\leq l$. Let $B$ be a maximum strong component of $D'$. We choose such a $T$ so that 

(1) $|B|$ is as large as possible. 

Let $P:=p_0p_1\cdots p_t$ ($t\geq2$) be a shortest dipath in $D$ such that $p_0,p_t\in B$ and $p_1,\cdots,p_{t-1}\in V(D)-B$ (by Lemma 4.1). We consider three cases in the following.

\noindent{\bf Case 1.} $t=2$.

By $t=2$, we have $p_1\in V(T)$. If $B=C_1$, then for any vertex $c_l\in C_l$, we have $|N_D^+(c_l)\backslash (B\cup p_1)|\geq m+1-1=m$ and $|N_D^+(c_l)\cap C_l|\geq m+1-|T|=1$. Let $(c_l,c_l')$ be an arc in $C_l$. Since $|N_D^+(c_l)\backslash (B\cup \{p_1,c_l'\})|\geq m+1-2=m-1$ and $|N_D^+(c_l')\backslash (B\cup \{p_1,c_l\})|\geq m+1-2=m-1$, by Lemma 4.3 (i),
we can find an out-double-star $T'\cong ODS(m;r,s)$ with center-arc $(c_l,c_l')$ such that $V(T')\cap (V(B)\cup V(P))=\O$. But then $V(B)\cup V(P)$ is contained in a strong component of $D-V(T')$, contrary to (1). If $B\neq C_1$, then for any vertex $c_1\in C_1$, we have $|N_D^-(c_1)\backslash (B\cup p_1)|\geq m+1-1=m$ and $|N_D^-(c_1)\cap C_1|\geq m+1-|T|=1$. Let $(c_1,c_1')$ be an arc in $C_1$. Since $|N_D^-(c_1)\backslash (B\cup \{p_1,c_1'\})|\geq m+1-2=m-1$ and $|N_D^-(c_1')\backslash (B\cup \{p_1,c_1\})|\geq m+1-2=m-1$, by Lemma 4.3 (ii), we can find an in-double-star $T''\cong IDS(m;r,s)$ with center-arc $(c_1,c_1')$ such that $V(T'')\cap (V(B)\cup V(P))=\O$. But then $V(B)\cup V(P)$ is contained in a strong component of $D-V(T'')$, contrary to (1).

\noindent{\bf Case 2.} $t=3$.

By the choice of $P$, we have $N_D^+(p_1)\cap B=\O$. Then $|N_D^+(p_1)\setminus (B\cup P)|\geq m+1-1=m$. Let $q\in N^+(p_1)\setminus (B\cup P)$.

\noindent{\bf Case 2.1.} $N_D^+(q)\cap B=\O$.

By $N_D^+(q)\cap B=\O$, we have $|N_D^+(q)\backslash (B\cup P)|\geq m+1-2=m-1$.
Let $w\in N_D^+(q)\setminus (B\cup P)$.

If $N_D^+(w)\cap B=\O$, then $|N_D^+(w)\backslash (B\cup P)|\geq m+1-2=m-1$.
Since $|N_D^+(q)\backslash (B\cup P \cup \{w\})|\geq m+1-3=m-2$ and $|N_D^+(w)\backslash (B\cup P \cup \{q\})|\geq m+1-3=m-2$, by Lemma 4.3 (i),
we can find an out-double-star $T'\cong ODS(m;r,s)$ with center-arc $(q,w)$ such that $V(T')\cap (V(B)\cup V(P))=\O$. But then $V(B)\cup V(P)$ is contained in a strong component of $D-V(T')$, contrary to (1).

If $N_D^+(w)\cap B\neq\O$, then $N_D^-(w)\cap B=\O$ (for otherwise, we can find a dipath $P'$ shorter that $P$). Thus $|N_D^-(w)\backslash (B\cup P)|\geq m+1-2=m-1$. Since $|N_D^+(q)\backslash (B\cup P \cup \{w\})|\geq m+1-3=m-2$ and $|N_D^-(w)\backslash (B\cup P \cup \{q\})|\geq m+1-3=m-2$, by Lemma 4.3 (iii), we can find an out-in-double-star $T'\cong OIDS(m;r,s)$ with center-arc $(q,w)$ such that $V(T')\cap (V(B)\cup V(P))=\O$. But then $V(B)\cup V(P)$ is contained in a strong component of $D-V(T')$, contrary to (1).

\noindent{\bf Case 2.2.} $N_D^+(q)\cap B\neq\O$.

By $N_D^+(q)\cap B\neq\O$, we have $N_D^-(q)\cap B=\O$ (for otherwise, we can find a dipath $P'$ shorter that $P$), and then $|N_D^-(q)\backslash (B\cup P)|\geq m+1-2=m-1$. Let $w'\in N_D^-(q)\setminus (B\cup P)$.

If $N_D^+(w')\cap B=\O$, then $|N_D^+(w')\backslash (B\cup P)|\geq m+1-2=m-1$.
Since $|N_D^-(q)\backslash (B\cup P \cup \{w'\})|\geq m+1-3=m-2$ and $|N_D^+(w')\backslash (B\cup P \cup \{q\})|\geq m+1-3=m-2$, by Lemma 4.3 (iii),
we can find an out-in-double-star $T'\cong OIDS(m;r,s)$ with center-arc $(w',q)$ such that $V(T')\cap (V(B)\cup V(P))=\O$. But then $V(B)\cup V(P)$ is contained in a strong component of $D-V(T')$, contrary to (1).

If $N_D^+(w')\cap B\neq\O$, then $N_D^-(w')\cap B=\O$ (for otherwise, we can find a dipath $P'$ shorter that $P$). Thus
$|N_D^-(w')\backslash (B\cup P)|\geq m+1-2=m-1$. Since $|N_D^-(q)\backslash (B\cup P \cup \{w'\})|\geq m+1-3=m-2$ and $|N_D^-(w')\backslash (B\cup P \cup \{q\})|\geq m+1-3=m-2$, by Lemma 4.3 (ii), we can find an in-double-star $T'\cong  IDS(m;r,s)$ with center-arc $(w',q)$ such that $V(T')\cap (V(B)\cup V(P))=\O$. But then $V(B)\cup V(P)$ is contained in a strong component of $D-V(T')$, contrary to (1).

\noindent{\bf Case 3.} $t\geq4$.

By the choice of $P$, we have $N_D^+(p_1)\cap (B\cup P)=\{p_2\}$. Then $|N_D^+(p_1)\setminus (B\cup P)|\geq m+1-1=m$. Let $q_1,\cdots,q_m\in N_D^+(p_1)\setminus (B\cup P)$. By the choice of $P$, we have $N_D^+(q_j)\cap (B\cup P)\subseteq\{p_1,p_2,p_3\}$ for each $j\in\{1,\cdots,m\}$. Thus $|N_D^+(q_j)\setminus(B\cup P)|\geq m+1-3=m-2$ for each $j\in\{1,\cdots,m\}$.
Let $w\in N_D^+(q_1)\setminus (B\cup P)$.

If $N_D^+(w)\cap B\neq\O$, then $N_D^-(w)\cap (B\cup P)\subseteq\{p_{t-2},p_{t-1}\}$ by the choice of $P$. Thus $|N_D^-(w)\setminus (B\cup P)|\geq m-1$. Since $|N_D^+(q_1)\setminus(B\cup P)|\geq m-2$ and $|N_D^-(w)\setminus (B\cup P)|\geq m-1$, by Lemma 4.3 (iii), we can find an out-in-double-star $T'\cong OIDS(m;r,s)$ with center-arc $(q_1,w)$ such that $V(T')\cap (V(B)\cup V(P))=\O$. But then $V(B)\cup V(P)$ is contained in a strong component of $D-V(T')$, contrary to (1).

If $N_D^-(w)\cap B\neq\O$, then $N_D^+(w)\cap (B\cup P)\subseteq\{p_1,p_2\}$ by the choice of $P$. Thus $|N_D^+(w)\setminus (B\cup P)|\geq m-1$. Since $|N_D^+(q_1)\setminus(B\cup P)|\geq m-2$ and $|N_D^+(w)\setminus (B\cup P)|\geq m-1$, by Lemma 4.3 (i), we can find an out-double-star $T'\cong ODS(m;r,s)$ with center-arc $(q_1,w)$ such that $V(T')\cap (V(B)\cup V(P))=\O$. But then $V(B)\cup V(P)$ is contained in a strong component of $D-V(T')$, contrary to (1).

In the following, we assume $N_D^+(w)\cap B=\O$ and $N_D^-(w)\cap B=\O$ for every $w\in N_D^+(q_1)\setminus (B\cup P)$. By the choice of $P$, we have $N_D^+(w)\cap (B\cup P)\subseteq\{p_1,p_2,p_3,p_4\}$ for each $w\in N_D^+(q_1)\setminus (B\cup P)$. Thus $|N_D^+(w)\setminus(B\cup P)|\geq m+1-4=m-3$ for each $w\in N_D^+(q_1)\setminus (B\cup P)$.

If there is some $j\in\{1,\cdots,m\}$, say $j=1$, such that $|N_D^+(q_1)\setminus (B\cup P)|\geq m-1$, let $w_1,\cdots,w_{m-1}\in N_D^+(q_1)\setminus (B\cup P)$. If there is some $k\in\{1,\cdots,m-1\}$, say $k=1$, such that $|N_D^+(w_1)\cap (B\cup P)|\geq m-2$, then by Lemma 4.3 (i), we can find an out-double-star $T'\cong ODS(m;r,s)$ with center-arc $(q_1,w_1)$ such that $V(T')\cap (V(B)\cup V(P))=\O$. But then $V(B)\cup V(P)$ is contained in a strong component of $D-V(T')$, contrary to (1). So we assume $|N_D^+(w_k)\setminus (B\cup P)|=m-3$ for each $k\in\{1,\cdots,m-1\}$. Then we have $N_D^+(w_k)\cap (B\cup P)=\{p_1,p_2,p_3,p_4\}$ for each $k\in\{1,\cdots,m-1\}$. By Lemma 4.3 (ii), there is an in-double-star $T''\cong IDS(m;r,s)$ with center-arc $(p_2,p_3)$ and vertex set $\{p_2,p_3,w_1,\cdots,w_{m-2}\}$ such that $V(T'')\cap (V(B)\cup V(P'))=\O$, where $P':=p_0p_1q_1w_{m-1}p_4\cdots p_t$. But then $V(B)\cup V(P')$ is contained in a strong component of $D-V(T'')$, contrary to (1).

Thus we assume that $|N_D^+(q_j)\setminus (B\cup P)|=m-2$ and $N_D^+(q_j)\cap (B\cup P)=\{p_1,p_2,p_3\}$ for each $j\in\{1,\cdots,m\}$. By $|N_D^+(q_j)\setminus (B\cup P)|=m-2$, there is an integer $j\in \{2,\cdots,m\}$, say $j=m$, such that $q_m\notin N_D^+(q_1)$. Since $\{q_2,\cdots,q_{m-1}\}\in N_D^-(p_2)$ and $|N_D^+(q_1)\setminus (B\cup P)|=m-2$, by Lemma 4.3 (iii), we can find an out-in-double-star $T'\cong OIDS(m;r,s)$ with center-arc $(q_1,p_2)$ such that $V(T')\cap (V(B)\cup V(P'))=\O$, where $P':=p_0p_1q_mp_3\cdots p_t$. But then $V(B)\cup V(P')$ is contained in a strong component of $D-V(T')$, contrary to (1).

The theorem is established as all cases lead to contradictions.  $\Box$

\end{document}